\makeatletter
\def\input@path{{D:/Stuff/Misc//}}
\makeatother
\documentclass[twocolumn,peerreview]{article}
\usepackage[T1]{fontenc}
\usepackage{verbatim}
\usepackage{float}
\usepackage{amsmath}
\usepackage{graphicx}

\makeatletter

\providecommand{\LyX}{L\kern-.1667em\lower.25em\hbox{Y}\kern-.125emX\@}

\onecolumn
\@addtoreset{equation}{section}

\makeatother
\begin{document}

\title{Dynamics of Jackson networks: perturbation theory}

\maketitle
Reuven Zeitak%
\footnote{reuven.zeitak@alcatel-lucent.co.il%
}

Alcatel Optical Networks Israel,2 Granit St., P.O.Box 7165,Petah Tikva,
Israel

\begin{abstract}
We introduce a new formalism for dealing with networks of queues.
The formalism is based on the Doi-Peliti second quantization method
for reaction diffusion systems. As a demonstration of the method's
utility we compute perturbatively the different time busy-busy correlations
between two servers in a Jackson network.
\end{abstract}
Jackson networks; second quantization; operator formalism; perturbation
theory;correlation functions;

\section{Introduction}

A Jackson network consists of a set of customers performing random
walks within a network of queues, where customers may also be added
(in a Poisson manner) or removed from the network. A well known theorem
(Jackson's theorem\cite{Jackson}) states that the instantaneous steady-state
behavior of such a queueing network has the same statistics as a set
of independent M/M/1 queues. In contrast to this simple result, dynamical
properties of queueing networks are quite involved due to the appearance
of correlations between events at different times. As very little
is known about the nature of these correlations, a systematic approach
would be interesting. Providing such an approach is one of the purposes
of this paper.

Another purpose of this paper is to point out a connection between
queueing networks and various physical models known as \char`\"{}reaction
diffusion models\char`\"{} (RD) that are used to model bulk chemical
reactions and various other particle dynamics. RD systems have benefited
enormously from a major insight into the problem that was made by
Doi and Peliti, who noticed ( independently) the usefulness of quantum
many body techniques to analyze the RD problem. Their insight was
that the RD problem could be written using quantum-mechanics-like
\char`\"{}second quantized\char`\"{} operators to describe the hopping
and interactions of the particles.

RD models are used to describe the microscopic motion of particles
through a medium which has a diffusive effect on the particles. If
one considers a queue customer as a particle and a server as a site
at which an interaction takes place, the hopping can be considered
as the result of the customers getting randomly routed to other servers
and being queued there. The description of a queueing network as a
set of customers performing random walks is more or less standard,
but the relationship to molecules adrift in a medium seems to have
been neglected. Once the analogy is made, we can borrow some of the
techniques used in these theories to reformulate some queueing models
in a suggestive (and in some cases simpler) form. This reformulation
will allow us to develop a systematic perturbation expansion of server
state correlations i.e. the busy-busy correlation of two different
queues (at different times). We shall calculate the first order term
in the perturbation, which appears to be a new result, and gives us
a formula that is found to be numerically valid across a wide range
of simulated network parameters.

The paper is organized as follows: after a brief review of relevant
previous work the Doi-Peliti formalism shall be introduced and applied
to queues in section \ref{sec:Second-quantization-formalism}. In
sub-section \ref{sub:Open-Jackson-network} we shall introduce the
dynamic operator that represents the Jackson network and show how
the equilibrium statistics can be easily derived from the operator.
The generalization to dynamic correlations will be discussed in section
\ref{sec:Queue-dynamics-in}. In section \ref{sec:Perturbation-expansion-for}
we shall generate a formal perturbation expansion for the \char`\"{}propagator\char`\"{}
of the Jackson network operator. We shall use the expansion to calculate
the (Laplace transformed) \char`\"{}busy-busy\char`\"{} correlation
function. In section \ref{sec:Comparisons-to-Simulations} the perturbation
results will be compared to simulations. Section \ref{sec:Conclusions}
is devoted to remarks on possible extensions and shortcomings of the
technique.

\section{\label{sec:Previous-Work}Previous Work}

Second quantization was introduced in reaction diffusion problems
by Doi and Peliti\cite{Doi+Peliti}. The method has been developed
quite extensively by Cardy and others\cite{Cardy}. The relationship
between RD and queueing was pointed out by\cite{RD+Queues}.

The possibility to compute correlation functions perturbatively for
the queue network depends on the knowledge of the Green's function
for the single M/M/1 queue. A convenient representation, that maps
into a normal-ordered second quantized form was derived by \cite{Greens}.

An approach that is very similar to the one presented here is that
of Massey\cite{Massey}\cite{Massey2}. In a series of papers he defines
\char`\"{}an operator theoretic approach\char`\"{} to Markovian queues.
These operators are reminiscent of the quantum mechanical (QM) creation
and annihilation operators introduced here. Indeed, the work presented
here could be viewed as an extension of Massey's work (although the
author came about the representation independently), but the emphasis
of this paper is quite different. This paper deals mostly with dynamical
properties of queues and the possibility to get perturbative results
for correlations. While section \ref{sec:Second-quantization-formalism}
can be viewed as a review of known queueing results in a second quantized
framework, the rest of the paper goes beyond this, to obtain new,
dynamical results. For readers unfamiliar with the operator theoretic
approach, section \ref{sec:Second-quantization-formalism} can be
used to make contact with standard queueing theory formalism.

\section{\label{sec:Second-quantization-formalism}Second quantization formalism}

\subsection{General}

We define the state of a queue as the number of particles (customers)
stored in the queue (we include the customer being serviced as being
the first in queue). A queue with $n$ particles is denoted by the
{}``ket'' $|n\rangle$. In addition, we define a set of orthogonal
{}``bra'' states $\langle m|$ such that the inner product$\langle m|n\rangle=\delta_{mn}$
. For a system of multiple queues indexed by $i$ a snapshot of the
system is given by a direct product of all single site states: $|n_{1},n_{2}\dots\rangle=\prod_{k}|n_{k}\rangle_{k}$.
The probability vector of all states can be given as a sum of all
possible configurations $|\psi\rangle=\sum_{n_{1}}\sum_{n_{1}}\dots P(n_{1},n_{2}\dots)\prod_{k}|n_{k}\rangle$.

Following the usual course\cite{Cardy}, we define creation and annihilation
operators (sometimes called {}``\emph{ladder}'' operators) $a^{+}$
and $a$ respectively, that have the following effect on the states:

\begin{eqnarray}
a|n\rangle &  & =n|n-1\rangle\,\,\,{\bf annihilation}\nonumber \\
a^{+}|n\rangle &  & =|n+1\rangle\,\,\,{\bf creation}\label{eq:General1}\end{eqnarray}

Which generalizes to 

\begin{eqnarray}
a_{i}|n_{1}\rangle_{1}|n_{2}\rangle_{2}\dots|n_{i}\rangle_{i}\dots= & n_{i}|n_{1}\rangle_{1}|n_{2}\rangle_{2}\dots|n_{i}-1\rangle_{i}\dots & \,\,\,{\bf annihilation}\nonumber \\
a_{i}^{+}|n_{1}\rangle_{1}|n_{2}\rangle_{2}\dots|n_{i}\rangle_{i}\dots= & |n_{1}\rangle_{1}|n_{2}\rangle_{2}\dots|n_{i}+1\rangle_{i}\dots & \,\,\,{\bf creation}\label{eq:general2}\end{eqnarray}
It turns out that although the annihilation operator $a$ has a simple
commutation relation with the creation operator $a^{+}$, it is not
a useful operator when dealing with the M/M/1 queue. This is due to
the fact that he probability to leave a queue is not a function of
the number of items in the queue (no mass action law\cite{MassAction}).
Instead of the $a$ operator, we shall define an modified annihilation
operator $Q$ such that 

\begin{eqnarray}
|n-1\rangle &  & =Q|n\rangle\,\,\,{\rm for}\,\, n>0\nonumber \\
0 &  & =Q|0\rangle\label{eq:Q for all N}\end{eqnarray}

So that $a^{+}$ can be considered as an operator that adds a single
client to a queue and $Q$ as an operator that removes a client (
i.e. via serving the client).

\subsection{\label{sub:Single-queue-M/M/1}Single queue M/M/1 at equilibrium}

Using the ladder operators $Q$ and $a^{+}$ \cite{qDeformed} , the
M/M/1 queue master equation 

\begin{eqnarray}
\dot{P}_{k} &  & =\mu(P_{k+1}-P_{k})+\lambda(P_{k-1}-P_{k})\,\,\, k>0\nonumber \\
\dot{P}_{0} &  & =\mu P_{1}-\lambda P_{0}\label{eq:}\end{eqnarray}

has the operator form of

\begin{alignat}{1}
\partial_{t}|\psi & \rangle={\cal L}|\psi\rangle\label{eq:MM1}\end{alignat}

with

\begin{equation}
{\cal L}=(1-a^{+})(\mu Q-\lambda)=\mu(1-a^{+})(Q-\rho).\label{eq:MM1operator}\end{equation}

Where $\rho=\frac{\lambda}{\mu}$ as usual.

There is nothing {}``mystical'' about the notation. The summation
of configurations that define the {}``wave function'' $|\psi\rangle$
is quite similar to a $z$-transform of the probability vector, i.e.
$\psi(\vec{z})=\sum_{n_{1}}\sum_{n_{1}}\dots P(n_{1},n_{2}\dots)\prod_{k}z_{k}^{n_{k}}$.
The ladder operator $a^{+}$ is equivalent to multiplying by $z$
and the lowering operator $a$ is equivalent to differentiating by
$z$ (that is $\partial_{z}$). The $Q$ operator is slightly less
familiar, but can be viewed as the operation of $\psi(z)\rightarrow\frac{\psi(z)-\psi(0)}{z}$.
It is just more convenient to treat these as abstract operators.

As pointed out in Massey's papers, the equilibrium solution 

\begin{alignat}{1}
0 & ={\cal L}|\psi\rangle\label{eq:MM1 Equilibrium}\end{alignat}

can be obtained by observing that in order to obtain $0=\mu(1-a^{+})(Q-\rho)|\psi\rangle$
it is enough to find a state that obeys $0=(Q-\rho)|\psi\rangle$.
We note that the modified annihilation operator $Q$ generates a set
of eigen-states, similar to the {}``coherent states'' that exist
in quantum mechanics. Namely

\begin{eqnarray}
x{\frac{1}{1-xa^{+}}}|0\rangle &  & =Q{\frac{1}{1-xa^{+}}}|0\rangle.\label{Qprop3}\end{eqnarray}

As can be seen by expanding in $x$. We see that $(Q-\rho){\frac{1}{1-xa^{+}}}|0\rangle=(x-\rho){\frac{1}{1-xa^{+}}}|0\rangle$
so that by setting $x=\rho$ we can solve equation \ref{eq:MM1 Equilibrium}
 with 

\begin{equation}
|\rho\rangle\equiv{\frac{1-\rho}{1-\rho a^{+}}}|0\rangle=\{1-\rho\}\{|0\rangle+\rho|1\rangle+\rho^{2}|2\rangle+\dots\}\label{eq:}\end{equation}

The formal time dependent solution to the M/M/1 queue, starting with
an initial probability state $|\psi_{0}\rangle$ is $|\psi(t)\rangle=e^{t{\cal L}}|\psi_{0}\rangle$
and due to the Markov-chain nature of the M/M/1 model we expect that
the long time behavior of arbitrary physical initial conditions relaxes
to the equilibrium state $|\rho\rangle$\cite{coherent}. In what
follows we shall sometimes call the time domain Green's function $e^{t{\cal L}}$
the propagator. Various average quantities can be obtained by constructing
expectation values with the {}``unit bra'' \begin{equation}
\langle I|\equiv\langle0|+\langle1|+\langle2|+...\label{eq:<I|}\end{equation}
 For example, suppose we are interested in the average busy ratio
$b$ of an equilibrized queue. The operator combination $a^{+}Q$
represents an object that returns $1$ when the queue is not empty
and $0$ otherwise. Hence $b=\langle I|a^{+}Q|\rho\rangle$ which
can be easily evaluated as $\rho$ .

\subsection{\label{sub:Open-Jackson-network}Open Jackson network of queues at
equilibrium}

The open Jackson network looks like this:

\begin{equation}
\partial_{t}|\psi\rangle=\sum_{j}(1-a_{j}^{+})[\sum_{i}(\delta_{ij}-r_{i\rightarrow j})\mu_{i}Q_{i}-\gamma_{j}]|\psi\rangle\label{eq:OpenJackson1}\end{equation}

To simplify, we define $L_{ij}=(\delta_{ij}-r_{i\rightarrow j})\mu_{i}$
and \begin{equation}
L_{ij}\rho_{i}=\gamma_{j}\label{eq:defz}\end{equation}

so that the open Jackson network operator becomes 

\begin{equation}
\partial_{t}|\psi\rangle=(1-a_{j}^{+})L_{ij}[Q_{i}-\rho_{i}]|\psi\rangle\label{eq:OpenJackson2}\end{equation}

In this form is now quite easy to verify that the product form

$|\rho_{1},\rho_{2},\dots\rangle=\prod_{k}\frac{(1-\rho_{k})}{1-a_{k}^{+}\rho_{k}}|0\rangle$
is the stationary solution of the open M/M/1 Jackson network, thus
proving Jackson's theorem.

\section{Queue dynamics in the second quantized formalism\label{sec:Queue-dynamics-in}}

\subsection{Expectation values and correlations}

We have seen that the second quantized formalism allows us to express
averages of functions of the occupation number as the expectation
values of various operators. For example, as stated above, the operator
$a^{+}Q$ projects out of a probability vector all states that are
non empty (namely $a^{+}Q|n\rangle=|n\rangle\textrm{ \ for }n>0$
and $a^{+}Q|0\rangle=0$). Similarly, $(a^{+})^{2}(Q)^{2}$ projects
out all states that have 2 or more customers in queue. In fact, the
term $(a^{+})^{n}Q^{n}-(a^{+})^{n+1}Q^{n+1}$projects out the state
that has exactly $n$ customers. We see that at least formally, given
a probability vector $|\psi\rangle$we can extract all queue occupation
information by considering expectation values of operators $O$, $\langle I|O|\psi\rangle$.

Suppose we know the queue state at time $0$ and we want to find the
average of some queue occupation number quantity at time $t$ . Assuming
that the initial queue state is $|\psi_{0}\rangle$, the formal solution
at time $t$ is $|\psi(t)\rangle=e^{t{\cal L}}|\psi_{0}\rangle$ so
that the expectation value of the operator $O$ becomes $\langle I|O|\psi(t)\rangle=\langle I|Oe^{t{\cal L}}|\psi_{0}\rangle$.
This expectation value can be interpreted as an average measurement
performed on a queue after it has evolved from the initial state for
a period $t$.

Let us now consider the case of correlation functions, i.e. measurements
taken at two different times. If we have two operators $O_{1},O_{2}$
that represent two measurements of the queue behavior, we can construct
the correlation function as $\langle I|O_{2}e^{tL}O_{1}|\psi\rangle$.
This object is equivalent to starting out with a distribution $|\psi\rangle$
of initial conditions of the queue, measuring the value of $O_{1}$,
then evolving the queue for a time $t$, measuring the value of $O_{2}$
and then averaging over all evolutions and all initial conditions
given by $|\psi\rangle$. If we are interested in steady state correlations,
we may replace the generic initial condition $|\psi\rangle$ with
the steady state solution $|\rho\rangle$ .

Correlation functions are important because they give dynamical information
regarding queue evolution. The correlation function is the first moment
of the joint probability distribution $<ab>=\int dadb\ abP(a,b)$,
so that the joint probability of our two measurements can be extracted
from a generating function $J(p,q)=\langle I|e^{ipO_{2}}e^{tL}e^{iqO_{2}}|\psi\rangle$.

\subsection{\label{sub:Single-queue-M/M/1}Single queue M/M/1 correlation functions}

What is the correlation between the server state (i.e. busy or empty)
at time $0$ and time $t$ ? If the initial queue state distribution
is $|\psi_{0}\rangle$ we need to compute

\begin{equation}
C(t)=\frac{\langle I|a^{+}Qe^{t{\cal L}}a^{+}Q|\psi_{0}\rangle}{\langle I|\psi_{0}\rangle}\label{eq:}\end{equation}

which simplifies for $|\psi_{0}\rangle=|\rho\rangle$ to $C(t)=\rho\langle I|Qe^{t{\cal L}}a^{+}|\rho\rangle$. 

Quantum mechanics teaches us that expressions such as this can be
dealt with conveniently if the operators inside are {}``normal-ordered'',
that is, brought to a form such that all annihilation operators are
near the ket $|\rho\rangle$ and all the creation operators are near
the bra $\langle I|$ . The propagator $g(t)=e^{t{\cal L}}$ does
not have a simple normal ordered form, but happily, its Laplace transform
$\hat{g}(\omega)=\frac{1}{\omega-{\cal L}}$ has been reduced to such
a form by \cite{Greens}, although they did not use an operator formalism.
For completeness we present the derivation in appendix \ref{sub:Single-queue-greens}.
The final result is

\begin{equation}
\hat{g}(\omega)=\frac{x/(\mu\rho)}{(1-xa^{+})}(1+\frac{x}{\rho-x}(1-a^{+}Q))\frac{1}{(1-xQ/\rho)}\label{eq:MM1 Greens}\end{equation}

where 

\begin{equation}
x(\omega)=\frac{(\omega/\mu+\rho+1)-\sqrt{(\omega/\mu+\rho+1)^{2}-4\rho}}{2}\label{eq:x(w)}\end{equation}
.

We note that $\hat{g}(\omega)$ is normal ordered and that computing
$\hat{C}(\omega)=\rho\langle I|Q\hat{g}(\omega)a^{+}|\rho\rangle$
becomes a simple (although tedious) exercise if one recalls that $Qa^{+}=1$.
The result is

\begin{eqnarray}
\rho\langle I|Q\hat{g}(\omega)a^{+}|\rho\rangle=\frac{x}{\mu(1-x)(\rho-x)}(\rho+\rho x-x)\label{eq:MM1Correlation}\end{eqnarray}

(see appendix \ref{sub:Single-queue-busy-busy} for details).

\subsection{\label{sub:Jackson-network-correlation}Jackson network correlation
functions}

Consider queue $\beta$ at time $0$ and queue $\alpha$ at time $t$
($\alpha\neq\beta$). If both queues are not empty at the respective
times, we say the busy-busy correlation is 1, and zero otherwise.
In equilibrium ( that is, starting out with the stationary state at
time zero), this function can be described as:

\begin{equation}
C_{\alpha\beta}(t)=\langle I|a_{\alpha}^{+}Q_{\alpha}e^{t{\cal L}}a_{\beta}^{+}Q_{\beta}|\rho\rangle=\rho_{\beta}\langle I|Q_{\alpha}e^{t{\cal L}}a_{\beta}^{+}|\rho\rangle\label{Cab}\end{equation}

where this time, ${\cal L}=(1-a_{j}^{+})L_{ij}[Q_{i}-\rho_{i}]$ is
the dynamic operator of the Jackson network. Unfortunately, the exact
normal ordered form for $e^{t{\cal L}}$ or for its Laplace transform
$\frac{1}{\omega-{\cal L}}$ is not known, so that in order to calculate
the correlation function we must resort to approximations. In the
next section we shall apply perturbation theory, which allows a systematic
approximation scheme in terms of a small parameter.

\section{\label{sec:Perturbation-expansion-for}Perturbation expansion for
Jackson networks}

\subsection{The General Formalism}

We shall consider a perturbative expansion around the diagonal part
of the Jackson operator. That is, we consider the Jackson network
as a perturbation around a set of independent M/M/1 queues, with a
weak coupling to each other. The source rate of each unperturbed queue
is $\mu_{k}\rho_{k}$, that is, the effective rate which appears in
the steady state solution to the full network. The perturbation is
the non Markovian effect of cross talk between the independent queues.

if we consider

\begin{equation}
L_{ij}=L^{0}+\epsilon L^{1}=(\delta_{ij})\mu_{i}-\epsilon r_{i\rightarrow j}\mu_{i}\label{eq:}\end{equation}

we may write the operator as 

\begin{eqnarray}
\begin{aligned}{\cal L}= & \sum_{i,j}(1-a_{j}^{+})(\delta_{ij}-\epsilon r_{i\rightarrow j})\mu_{i}[Q_{i}-\rho_{i}]=\\
= & \sum_{j}(1-a_{j}^{+})\mu_{j}[Q_{j}-\rho_{j}]-\epsilon\sum_{i,j}(1-a_{j}^{+})r_{i\rightarrow j}\mu_{i}[Q_{i}-\rho_{j}]=\\
= & L^{0}+\epsilon W\end{aligned}
\label{eq:}\end{eqnarray}

Where we have set $L^{0}=\sum_{j}(1-a_{j}^{+})\mu_{j}[Q_{j}-\rho_{j}]$
as the unperturbed operator and $\textrm{W}=-\sum_{i,j}(1-a_{j}^{+})r_{i\rightarrow j}\mu_{i}[Q_{i}-\rho_{j}]$
the perturbation.

While it may seem unnatural to fix the values of the $\rho_{k}$ and
to expand around the diagonal terms of $L$ rather than fixing the
values of the $\gamma_{k}$, this form of splitting of the operator
has the advantage that the stationary solution remains unchanged for
all values of $\epsilon$. When we calculate correlations with respect
to the stationary state, we can avoid modifications (familiar in perturbation
theoretic expansions) due to the change in the stationary state.

The propagator for the unperturbed operator is just a product of propagators
for uncoupled queues $\mathcal{G}^{0}(t)=\prod g_{i}^{0}(t)$

The propagator $e^{t{\cal L}}$ becomes 

\begin{equation}
e^{t[L^{0}+\epsilon W]}=\mathcal{G}^{0}(t)+\epsilon\int_{0}^{t}d\tau\mathcal{G}^{0}(\tau)W\mathcal{G}^{0}(t-\tau)+\epsilon^{2}\int_{0}^{t}d\tau_{1}\int_{0}^{\tau_{1}}d\tau_{2}\mathcal{G}^{0}(\tau_{1})W\mathcal{G}^{0}(\tau_{2})W\mathcal{G}^{0}(t-\tau_{1}-\tau_{2})+...\label{eq:PropPT_t}\end{equation}

and Laplace transforming yields:

\begin{equation}
G(\omega)=G^{0}(\omega)+\epsilon G^{0}(\omega)WG^{0}(\omega)+\epsilon^{2}G^{0}(\omega)WG^{0}(\omega)WG^{0}(\omega)+...\label{eq:PropPT_w}\end{equation}

While trying to obtain useful perturbative results we are faced with
two technical difficulties, the first being the need to put the terms
in the perturbation series in normal order and the second is the need
to Laplace transform a product of propagators. The transformation
of the product form of $\mathcal{G}^{0}(t)=\prod g_{i}^{0}(t)$ into
the convolved Laplace transformed form of $G^{0}(\omega)$ is rather
involved. Dealing only with off-diagonal correlations and only the
first order in the perturbation expansion will allow us to simplify
matters, as long as we consider direct calculations of correlation
functions.

\subsection{Perturbation expansion of the busy-busy correlation function}

Consider queue $\beta$ at time $0$ and queue $\alpha$ at time $t$
($\alpha\neq\beta$). If both queues are not empty at the respective
times, we say the busy-busy correlation is 1, and zero otherwise.
As pointed out in \ref{sub:Jackson-network-correlation} the busy
busy correlation function is:

\begin{equation}
C_{\alpha\beta}(t)=\langle I|a_{\alpha}^{+}Q_{\alpha}\mathcal{G}(t)a_{\beta}^{+}Q_{\beta}|\rho\rangle=\rho_{\beta}\langle I_{\alpha}|Q_{\alpha}\mathcal{G}(t)a_{\beta}^{+}|\rho_{\alpha}\rangle\label{eq:}\end{equation}

Expanding $C_{\alpha,\beta}(t)$in $\epsilon$ as $C_{\alpha,\beta}(t)=C_{\alpha,\beta}^{0}(t)+\epsilon C_{\alpha,\beta}^{1}(t)+\dots$
we find (using the propagator expansion in the time domain (equation
(\ref{eq:PropPT_t})),

\begin{eqnarray}
C_{\alpha,\beta}^{0}(t) & = & \rho_{\beta}\langle I_{\alpha}|Q_{\alpha}\mathcal{G}^{0}(t)a_{\alpha}^{+}|\rho_{\alpha}\rangle\nonumber \\
C_{\alpha,\beta}^{1}(t) & = & \rho_{\beta}\int_{0}^{t}\langle I|Q_{\alpha}\mathcal{G}^{0}(\tau)W\mathcal{G}^{0}(t-\tau)a_{\beta}^{+}|\rho\rangle\label{eq:}\end{eqnarray}

Furthermore, the product structures of ${\cal G}^{0}$ and of the
stationary state $|\rho\rangle$ allow us to decouple the zero order
term in the expansion:

\begin{equation}
C_{\alpha,\beta}^{0}(t)=\rho_{\beta}\langle I|Q_{\alpha}g_{\alpha}^{0}(t)g_{\beta}^{0}(t)a_{\beta}^{+}|\rho\rangle=\rho_{\beta}\langle I|Q_{\alpha}a_{\beta}^{+}|\rho\rangle=\rho_{\beta}<Q_{\alpha}><a_{\beta}^{+}>\label{eq:}\end{equation}

and the first order perturbation is:

\begin{equation}
C_{\alpha,\beta}^{1}(t)=\rho_{\beta}\sum_{ij}\delta L_{i,j}\int_{0}^{t}\langle I|Q_{\alpha}g_{\alpha}^{0}(\tau)(1-a_{j}^{+})[Q_{i}-\rho_{i}]g_{\beta}^{0}(t-\tau)a_{\beta}^{+}|\rho\rangle d\tau.\label{eq:}\end{equation}

All cases where $\alpha\neq j$or $\beta\neq i$ vanish due to the
commutations $[X_{i},Y_{j}]=0$ for $i\neq j$, leaving us with

\begin{equation}
C_{\alpha,\beta}^{1}(t)=\rho_{\beta}\delta L_{\beta,\alpha}\int_{0}^{t}\langle I|Q_{\alpha}g_{\alpha}^{0}(\tau)(1-a_{\alpha}^{+})[Q_{\beta}-\rho_{\beta}]g_{\beta}^{0}(t-\tau)a_{\beta}^{+}|\rho\rangle d\tau.\label{eq:}\end{equation}

Again, the product structure of the stationary solution comes to our
aid allowing us to decouple:

\[
C_{\alpha,\beta}^{1}(t)=\rho_{\beta}\delta L_{\beta,\alpha}\int_{0}^{t}\langle I_{\alpha}|Q_{\alpha}g_{\alpha}^{0}(\tau)(1-a_{\alpha}^{+})|\rho_{\alpha}\rangle\langle I_{\beta}|[Q_{\beta}-\rho_{\beta}]g_{\beta}^{0}(t-\tau)a_{\beta}^{+}|\rho_{\beta}\rangle d\tau.\]

This is simple enough to allow us to perform a Laplace transform:

\[
\hat{C}{}_{\alpha\beta}(\omega)=\hat{C}{}_{\alpha,\beta}^{0}(\omega)+\epsilon\hat{C}{}_{\alpha,\beta}^{1}(\omega)+\dots.\]

with

\begin{eqnarray*}
\hat{C}{}_{\alpha,\beta}^{0}(\omega)= & \frac{\rho_{\beta}\langle I|Q_{\alpha}a_{\beta}^{+}|\rho\rangle}{\omega}=\frac{\rho_{\alpha}\rho_{\beta}}{\omega}\\
\hat{C}{}_{\alpha,\beta}^{1}(\omega)= & \rho_{\beta}\delta L_{\beta,\alpha}\langle I_{\alpha}|Q_{\alpha}\hat{g}{}_{\alpha}^{0}(\omega)(1-a_{\alpha}^{+})|\rho_{\alpha}\rangle\langle I_{\beta}|[Q_{\beta}-\rho_{\beta}]\hat{g}{}_{\beta}^{0}(\omega)a_{\beta}^{+}|\rho_{\beta}\rangle\end{eqnarray*}

Expanding the $\hat{C}{}_{\alpha,\beta}^{1}(\omega)$ expression we
remain with

\[
\hat{C}{}_{\alpha\beta}(\omega)=\frac{\rho_{\alpha}\rho_{\beta}}{\omega}+\epsilon\rho_{\beta}\delta L_{\beta,\alpha}[\frac{\rho_{\alpha}}{\omega}-\langle I_{\alpha}|Q_{\alpha}\hat{g}{}_{\alpha}^{0}(\omega)a_{\alpha}^{+}|\rho_{\alpha}\rangle][\langle I_{\beta}|Q_{\beta}\hat{g}{}_{\beta}^{0}(\omega)a_{\beta}^{+}|\rho_{\beta}\rangle-\frac{\rho_{\beta}}{\omega}]...\]

and using the results in appendix \ref{sub:Single-queue-busy-busy}

\begin{equation}
\hat{C}{}_{\alpha\beta}(\omega)=\frac{\rho_{\alpha}\rho_{\beta}}{\omega}+\epsilon\rho_{\beta}\delta L_{\beta,\alpha}[\frac{\rho_{\alpha}-x_{\alpha}}{\omega}-\frac{x_{\alpha}}{\mu_{\alpha}\rho_{\alpha}(1-x_{\alpha})}][\frac{x_{\beta}}{\mu_{\beta}\rho_{\beta}(1-x_{\beta})}+\frac{x_{\beta}-\rho_{\beta}}{\omega}]+...\label{eq:MainResult}\end{equation}

This is the main result of the paper. The non zero $O(\epsilon)$
term in the expansion is a direct demonstration of the fact that the
actual dynamics of Jackson networks are not equivalent to an independent
Poissonian arrival system, even though the stationary state appears
to behave as one. In the next section we shall compare our result
to simulations of various networks. But before that we shall comment
on the correlation behavior a single queue with respect to itself.

\subsection{Same queue correlations and busy periods }

The same queue correlation function 

$C_{\alpha\alpha}(t)=\langle I|a_{\alpha}^{+}Q_{\alpha}\mathcal{G}(t)a_{\alpha}^{+}Q_{\alpha}|\rho\rangle=\rho_{\alpha}\langle I|Q_{\alpha}\mathcal{G}(t)a_{\alpha}^{+}|\rho\rangle$
cannot be deduced as special case of the inter-queue correlation calculated
above. This is due to two facts. The first is that the unperturbed
situation $C_{\alpha\alpha}^{0}(t)$ is different: $C_{\alpha\alpha}^{0}(t)=\rho_{\alpha}\langle I|Q_{\alpha}g_{\alpha}^{0}(t)a_{\alpha}^{+}|\rho\rangle$
, depending on one $g^{0}$ only. The second fact is that the first
order term in $\epsilon$ vanishes because the perturbation term $W$
doesn't contribute to the diagonal, $W_{\alpha\alpha}=0$.

It is beyond the scope of this paper to deal with the second order
contribution, however, an interesting result can be extracted from
this behavior- the mean busy period of a single server in a Jackson
network is identical to the independent server mean busy period. This
is deduced by considering the limit of two nearby (in time) measurements:
The busy state of a server at time $t$ and the busy state of the
server at time $t+\tau$ where $\tau$ is small. If the first measurement
is {}``busy'' and the second is {}``idle'', a busy period has
terminated somewhere between $t$ and $t+\tau$. In the limit of $\tau\rightarrow0$
this will give $\tau$ times the density of busy period endings. On
the other hand, the situation in which these to events occur is exactly
what is measured by the correlation

$\langle I|(1-a_{\alpha}^{+}Q_{\alpha})\mathcal{G}(\tau)a_{\alpha}^{+}Q_{\alpha}|\rho\rangle=\langle I|(1-a_{\alpha}^{+}Q_{\alpha})\left[1+\tau L\right]a_{\alpha}^{+}Q_{\alpha}|\rho\rangle+O(\tau^{2})$

However, $L=L^{0}+\epsilon W$ and because the first term in $\epsilon$
vanishes, we are left with the result,

$\langle I|(1-a_{\alpha}^{+}Q_{\alpha})\mathcal{G}(\tau)a_{\alpha}^{+}Q_{\alpha}|\rho\rangle=\tau\langle I|(1-a_{\alpha}^{+}Q_{\alpha})L^{0}a_{\alpha}^{+}Q_{\alpha}|\rho\rangle+O(\tau^{2})$,
exactly the same result as for the unperturbed case.

The density of the busy period endings is the inverse of the mean
time between two busy periods, which is composed of the sum of the
mean busy period and the mean idle period. Since the fraction of the
of time that the queue is busy is also independent of the perturbation,
the mean busy period must remain unchanged, and this result holds
for all orders of the perturbation.

\section{Comparisons to Simulations\label{sec:Comparisons-to-Simulations}}

\subsection{On the fly Laplace transforms}

In order to test the accuracy of the perturbation expansion we have
used simulations to compute busy-busy correlations for various Jackson
networks. We simulate a Jackson network and monitor the queue states.
We perform the Laplace transforms on the fly by considering the product
of the instantaneous server state of queue $\alpha$ with the exponential
averaged server state of queue $\beta$. The cumulative average of
these products are equivalent to the Laplace transform of the correlation
functions at equilibrium. The simulation used is event driven so that
it is useful to describe how these averages can be performed during
an event driven simulation.

If we term the {}``busy'' indicator of the server of queue $\eta$
at time $t$ by $b_{\eta}$ , We consider

\begin{equation}
B_{\beta}(\omega,T)=\int_{0}^{\infty}b_{\beta}(T-t)e^{-\omega t}dt\label{eq:OntheFly1}\end{equation}

and

\begin{flushleft}
\begin{equation}
C_{\alpha\beta}(\omega,T)=\int_{0}^{\infty}b_{\alpha}(T)b_{\beta}(T-t)e^{-\omega t}dt=b_{\alpha}(T)B_{\beta}(\omega,T)\label{eq:OnTheFly1.1}\end{equation}

\par\end{flushleft}

\begin{flushleft}
so that $C_{\alpha\beta}(\omega)=<C_{\alpha\beta}(\omega,T)>_{T}=<b_{\alpha}(T)B_{\beta}(\omega,T)>_{T}$.
\par\end{flushleft}

So that $B_{\beta}(T)$can be obtained by the exponential averaging
defined by

\begin{equation}
B(T_{0}+\Delta)=e^{-\omega\Delta}[B(T_{0})+\int_{0}^{\Delta}e^{\omega\tau}b(T_{0}+\tau)d\tau]\label{eq:OnTheFly3}\end{equation}

Furthermore, the busy status of the servers remain fixed between each
packet arrival (or departure) event, thus the exponential average
is very easy to evaluate in an event driven simulation. Suppose that
there is no change in the busy state of either queue between times
$T_{0}$ and $T_{1}$. Then,

\begin{equation}
B_{\beta}(T_{1})=e^{-\omega(T_{1}-T_{0})}B(T_{0})+\frac{(1-e^{-\omega(T_{1}-T_{0})})}{\omega}b_{\beta}(T_{0}+)\label{eq:OnTheFly4}\end{equation}

Averages can be obtained by integrating $C_{\alpha\beta}(\omega,T)$
over time and normalizing. Consider the contribution to the integral
of $C_{\alpha\beta}(\omega,T)$ between the times $T_{0}$ and $T_{1}=T_{0}+\Delta$:

\begin{equation}
\int_{0}^{\Delta}C_{\alpha\beta}(\omega,T_{0}+\tau)d\tau=\int_{0}^{\Delta}b_{\alpha}(T_{0}+\tau)B_{\beta}(\omega,T_{0}+\tau)d\tau\label{eq:OnTheFly5}\end{equation}

but since $b_{\alpha},b_{\beta}$ are fixed during this period,

\[
\int_{0}^{\Delta}b_{\alpha}(T_{0}+\tau)B_{\beta}(\omega,T_{0}+\tau)d\tau=b_{\alpha}(T_{0}+)\int_{0}^{\Delta}[e^{-\omega\tau}B_{\beta}(T_{0}+)+\frac{(1-e^{-\omega\tau})}{\omega}b_{\beta}(T_{0}+)]d\tau=\]

\[
=b_{\alpha}(T_{0}+)B_{\beta}(T_{0}+)\frac{(1-e^{-\omega\Delta})}{\omega}+\frac{b_{\alpha}(T_{0}+)b_{\beta}(T_{0}+)}{\omega}[\Delta-\frac{(1-e^{-\omega\Delta})}{\omega}]\]

This gives us two equations to update at each event:

\begin{eqnarray}
B_{\beta}(T_{1}) & = & e^{-\omega\Delta}B(T_{0})+\frac{(1-e^{-\omega\Delta})}{\omega}b_{\beta}(T_{0}+\epsilon)\nonumber \\
\int_{0}^{T_{1}}dtC_{\alpha\beta}(t) & = & \int_{0}^{T_{0}}dtC_{\alpha\beta}(t)+b_{\alpha}(T_{0}+)B_{\beta}(T_{0}+)\frac{(1-e^{-\omega\Delta})}{\omega}+\nonumber \\
 &  & \frac{b_{\alpha}(T_{0}+)b_{\beta}(T_{0}+)}{\omega}[\Delta-\frac{(1-e^{-\omega\Delta})}{\omega}]\label{eq:OnTheFly6}\end{eqnarray}

(the $\omega$ dependence is suppressed for brevity).

Normalizing to $\frac{1}{T}\int_{0}^{T}dtC_{\alpha\beta}(\omega,t)$
gives an estimate of $<C_{\alpha\beta}(\omega,T)>_{T}$ and thus of
$C_{\alpha\beta}(\omega)$.

\subsection{Simulation results}

\begin{figure}[H]
\includegraphics[scale=1.6]{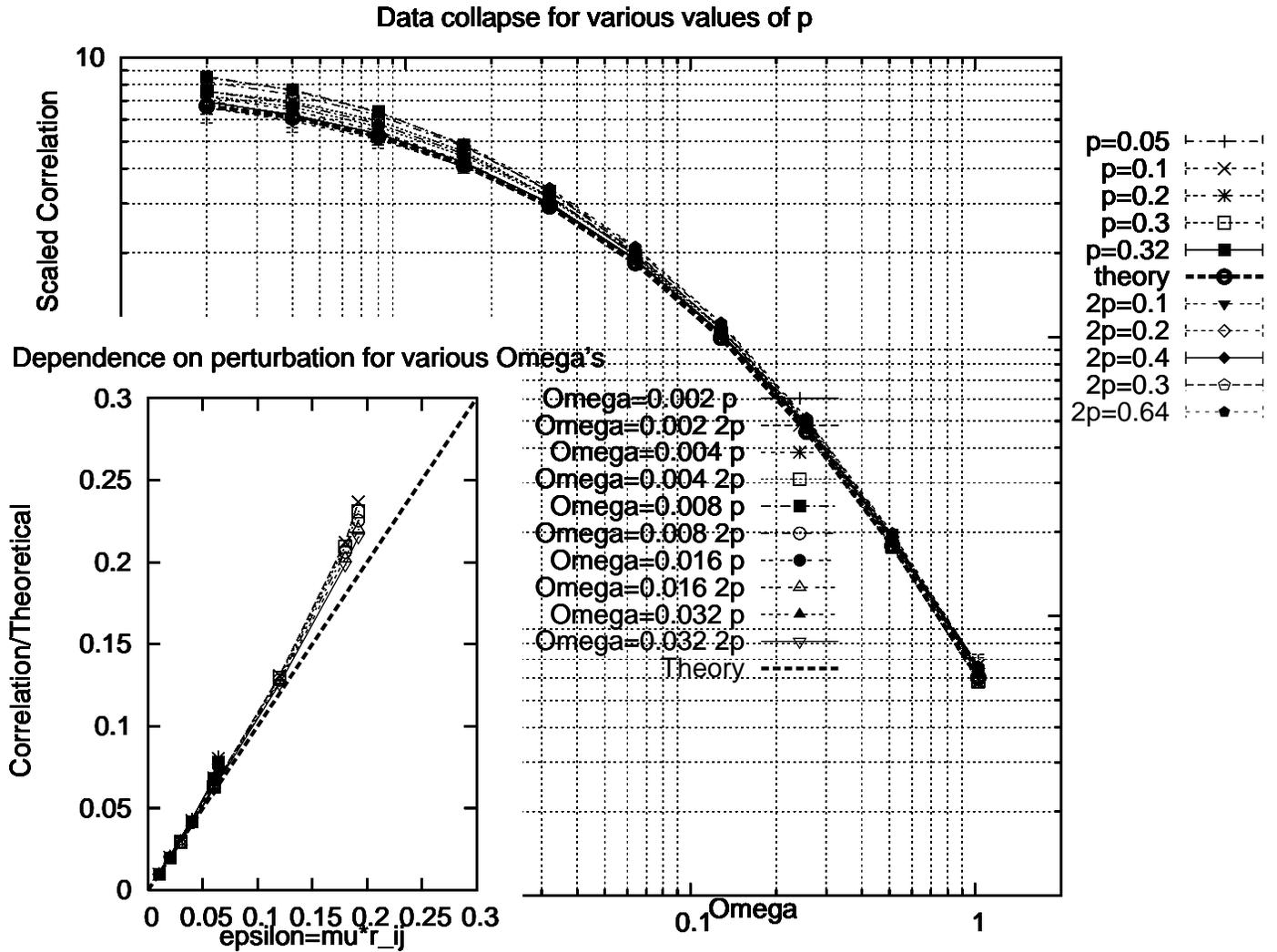}

\caption{\noindent \label{fig:Data-collapse-for}Type (1) data collapse for
the Laplace transform of the (subtracted) correlation function, scaled
by the perturbation value. the lines marked by 2p represent $<b{}_{2}b_{1}>$
correlations, and the lines marked by p represent $<b_{1}b_{2}>$
correlations. The deviations from the first order result (marked {}``theory'')
are due to higher order contributions. Error bars mark one RMS of
statistical deviation ( based on splitting the simulation run into
sub runs of 10,000 seconds). The total simulation is approximately
2,560,000 seconds , or 1M events on queue 1 and about 1.7M events
on queue 2. The inset shows type (2) data collapse- dividing out the
(unscaled) perturbation result and plotting the result versus the
perturbation amplitude. A $45^{o}$ straight line (marked {}``Theory'')
indicates exact matching to first order theory. Error bars are suppressed
for clarity.}

\end{figure}

We simulated a set of networks, for an arbitrary choice of $\rho$
(defined by equation (\ref{eq:defz}) ). we measured off-diagonal
correlations for various values of $\omega$ and subtracted the $0^{th}$
order perturbation term $\frac{\rho_{\alpha}\rho_{\beta}}{\omega}$
from the measured values. The data obtained is presented in two ways:

\begin{enumerate}
\item For networks that only differ by the value of the perturbation, we
plot the value of the subtracted correlation, normalized by the perturbation
$\delta L_{\beta,\alpha}=r_{\beta\rightarrow\alpha}\mu_{\beta}$ ,
as a function of $\omega$. We expect that the data collapse to the
same curve as long as first order perturbation theory is accurate.
\item In order to allow for data collapse of different families of networks,
the subtracted values were normalized by the computed perturbation
(equation \ref{eq:MainResult}) $\rho_{\beta}[\frac{\rho_{\alpha}-x_{\alpha}}{\omega}-\frac{x_{\alpha}}{\mu_{\alpha}\rho_{\alpha}(1-x_{\alpha})}][\frac{x_{\beta'}}{\mu_{\beta}\rho_{\beta}(1-x_{\beta})}+\frac{x_{\beta}-\rho_{\beta}}{\omega}]$,
leaving us with values were plotted with respect to $\delta L_{\beta,\alpha}=r_{\beta\rightarrow\alpha}\mu_{\beta}$
. Accurate results are reflected by a straight line collapse.
\end{enumerate}
A family of 2X2 matrices was tested: $r_{i\rightarrow j}=\left[\begin{array}{cc}
0 & p\\
2p & 0\end{array}\right]$, fixing $\{\rho_{1},\rho_{2}\}=\{0.3,0.7\}$ and $\{\mu_{1},\mu_{2}\}=\{0.3,0.2\}$,
we scanned the range with $p\in[0.05,0.1,0.2,0.3,0.32]$ and plotted
the results ( scaled by the perturbation) in figure \ref{fig:Data-collapse-for}. 

While a good fit is observed for all values of the perturbation and
of the frequency, systematic deviations, probably due to higher order
corrections are observable. Fits of similar quality were obtained
for larger networks (data not shown).

\section{\label{sec:Conclusions}Conclusions}

Having described the utility of the second quantization point of view
it might be considered natural to point out other queueing quantities
that could be calculated by various manipulations of the technique;
For example, different correlation functions such as the queue depth
correlator and other moments of the queue depth can be obtained. Similarly,
obtaining higher order terms in the perturbation expansion appears
to be a technically non trivial extension. Indeed, a diagrammatic
expansion for managing higher order terms may be formulated\cite{diagrams}.
The dynamics of other queueing models besides M/M/1 also seem to be
approachable in this technique (as pointed out above, the M/M/$\infty$
queue and network can be written using QM harmonic-oscillator ladder
operators). However, instead of dwelling on various extensions of
the formalism (which will be explored elsewhere), the author would
like to point out some shortcomings of the technique, in other words-
what is lacking in the second quantization formalism. 

The reason a {}``reaction diffusion'' approach is useful for queues
is due to the fact the the questions asked in this paper do not relate
directly to the {}``FIFO'' nature of the problem. The particles
were viewed as indistinguishable, and the entire information about
a queue state was encoded in the number of particles waiting at the
queue. The information about the $\textit{order}$ of the particles
waiting in the queue is lost in this representation. Thus, questions
that relate to the experience of a specific particle (e.g. waiting
times) cannot be formulated. Issues such as {}``jitter'' ( how the
inter-arrival time of two particles at the destination is related
to the {}``inter-injection'' time at the source) cannot be addressed
either.

This is somewhat analogous to the difference between {}``Eulerian''
and {}``Lagrangian'' views in fluid dynamics\cite{tritton}. In
the Lagrangian approach one considers the trajectories of tagged fluid
parcels that are advected with the flow. The Eulerian approach examines
the behavior at fixed positions in space and considers the density
and velocity of the fluid flowing through these positions. It could
be argued that while the Eulerian approach is technically more understood,
the Lagrangian approach captures various aspects that are very difficult
to introduce in the Eulerian point of view (e.g. the invariance under
Galilean invariance that removes the {}``sweeping'' effect that
masks various correlation functions).

It is difficult to ask {}``Lagrangian'' type questions it the operator
formalism. Massey attempted an extension in \cite{Massey2}, but in
is unclear how to use the results obtained there. In the framework
of the asymmetric simple exclusion process(ASEP) it is possible to
examine the behavior of {}``tracer'' particles\cite{ASEP}, but
the application to queues is not obvious.

\appendix

\section{Single queue green's function\label{sub:Single-queue-greens}}

The propagator is defined as the operator $g(t)=e^{t{\cal L}}$, and
the Green's function is its Laplace transform $\hat{g}(\omega)=\frac{1}{\omega-{\cal L}}$.
We are interested in a {}``normal ordered form'' for the Green's
function in which all the $Q$ operators are on the right and the
$a^{+}$are on the left ( so that $g|0\rangle$ is a simple calculation).
Following \cite{Greens}, we start with 

\begin{equation}
\hat{g}(\omega)=\frac{1}{\omega-\mu(1-a^{+})(Q-z)}=\frac{1}{\omega-\mu(Q-\rho+a^{+}\rho-1)+\mu(a^{+}Q-1)}.\label{eq:Inv1}\end{equation}

The authors of \cite{Greens} note that $\hat{g_{0}}(\omega)=\frac{1}{\omega-\mu(Q-\rho+a^{+}\rho-1)}$
is easily normal ordered by the ansatz: $\frac{x/(\mu\rho)}{(1-xa^{+})}\frac{1}{(1-xQ/\rho)}$
(for a yet to be determined value of $x$), which can be seen by directly
applying

$[\omega-\mu(Q-\rho+a^{+}\rho-1)]\frac{1}{(1-xa^{+})}\frac{1}{(1-xQ/\rho)}=$

$=[\omega-\mu(x-\rho-1)-\frac{\mu\rho}{x}]\frac{1}{(1-xa^{+})}\frac{1}{(1-xQ/\rho)}+\frac{\mu\rho}{x}[1-xQ/\rho]\frac{1}{(1-xQ/\rho)}=$

\begin{equation}
[\omega-\mu(x-\rho-1)-\frac{\mu\rho}{x}]\frac{1}{(1-xa^{+})}\frac{1}{(1-xQ/\rho)}+\frac{\mu\rho}{x}.\label{eq:G0}\end{equation}

If we demand that $x(\omega)$solves

\begin{equation}
[\omega-\mu(x-\rho-1)-\frac{\mu\rho}{x}]=0,\label{eq:ForX}\end{equation}

then the above inversion implies that we can write $\hat{g_{0}}(\omega)$
as

\begin{equation}
\hat{g_{0}}(\omega)=\frac{x/(\mu\rho)}{(1-xa^{+})}\frac{1}{(1-xQ/\rho)}.\label{eq:G0normal}\end{equation}

Given $\hat{g_{0}}(\omega)$ the authors of \cite{Greens} perturbatively
iterate for $G$:\begin{eqnarray*}
\hat{g}(\omega) & = & \frac{1}{\hat{g}_{0}^{-1}+\mu(a^{+}Q-1)}=\\
= &  & \hat{g}_{0}-\hat{g}_{0}\mu(a^{+}Q-1)\hat{g}_{0}+\hat{g}_{0}\mu(a^{+}Q-1)\hat{g}_{0}\mu(a^{+}Q-1)\hat{g}_{0}+...\,\end{eqnarray*}

noting that

\[
(a^{+}Q-1)\hat{g}_{0}(a^{+}Q-1)=|0\rangle\langle0|\frac{x/(\mu\rho)}{(1-xa^{+})}\frac{1}{(1-xQ/\rho)}|0\rangle\langle0|=(1-a^{+}Q)x/(\mu z)\]

we retrieve

$\hat{g}(\omega)=\hat{g}_{0}-\frac{\mu\rho}{\rho-x}\hat{g}_{0}(a^{+}Q-1)\hat{g}_{0}=\hat{g}_{0}-\frac{x^{2}/(\mu\rho)}{\rho-x}\frac{}{(1-xa^{+})}(a^{+}Q-1)\frac{1}{(1-xQ/\rho)}=$

$=\frac{x/(\mu\rho)}{(1-xa^{+})}\frac{1}{(1-xQ/\rho)}-\frac{x^{2}/(\mu\rho)}{\rho-x}\frac{}{(1-xa^{+})}(a^{+}Q-1)\frac{1}{(1-xQ/\rho)}=$

\begin{equation}
=\frac{x/(\mu\rho)}{(1-xa^{+})}(1-\frac{x}{\rho-x}(a^{+}Q-1))\frac{1}{(1-xQ/\rho)}.\label{eq:Gnormal}\end{equation}

We are left with solving equation \ref{eq:ForX}  for $x(\omega)$:

$x_{\pm}=\frac{(\omega/\mu+\rho+1)\pm\sqrt{(\omega/\mu+\rho+1)^{2}-4\rho}}{2}$

in order to decide the sign of the root, we note that for convergence
we must have $|x|<1$. We further note that for $\omega=0$

$x_{\pm}=\frac{(\rho+1)\pm\sqrt{(\rho+1)^{2}-4\rho}}{2}=\frac{(\rho+1)\pm|\rho-1|}{2}=\{\rho,1\}$
$(\rho<1)$

and since we expect convergence at this point, we must choose the
negative root

\begin{equation}
x(\omega)=\frac{(\omega/\mu+\rho+1)-\sqrt{(\omega/\mu+\rho+1)^{2}-4\rho}}{2}\label{eq:x(w)}\end{equation}

\section{\label{sub:Single-queue-busy-busy}Single queue busy-busy correlation}

We want to evaluate $\langle I|Q\hat{g}(\omega)a^{+}|\rho\rangle$.
The idea is to push all the $Q$'s to the right hand side and all
the $a^{+}$'s to the left.

\begin{equation}
\langle I|Q\hat{g}(\omega)a^{+}|\rho\rangle=\frac{x}{\mu\rho}\langle I|Q\frac{}{(1-xa^{+})}(1-\frac{x}{\rho-x}(a^{+}Q-1))\frac{1}{(1-xQ/\rho)}a^{+}|\rho\rangle\label{eq:busybusy}\end{equation}

$=\frac{x}{\mu\rho}\langle I|[Q+\frac{x}{(1-xa^{+})}](1-\frac{x}{\rho-x}(a^{+}Q-1))[\frac{x/\rho}{(1-xQ/\rho)}+a^{+}]|\rho\rangle$

$=\frac{x}{\mu\rho}\langle I|[Q+\frac{x}{(1-x)}](1-\frac{x}{\rho-x}(a^{+}Q-1))[\frac{x/\rho}{(1-x)}+a^{+}]|\rho\rangle$

$=\frac{x}{\mu\rho}\langle I|[Q+\frac{x}{(1-x)}]([\frac{x/\rho}{(1-x)}+a^{+}]-\frac{x}{\rho-x}(a^{+}Q-1)[\frac{x/\rho}{(1-x)}+a^{+}])|\rho\rangle$

$=\frac{x}{\mu\rho}\langle I|[Q+\frac{x}{(1-x)}]([\frac{x/\rho}{(1-x)}+a^{+}]-\frac{x}{\rho-x}[(a^{+}\rho-1)\frac{x/\rho}{(1-x)}])|\rho\rangle$

$=\frac{x}{\mu\rho}\langle I|[Q+\frac{x}{(1-x)}]([\frac{x}{(1-x)(\rho-x)}+a^{+}\frac{(\rho-\rho x-x)}{(\rho-x)(1-x)})|\rho\rangle$

$=\frac{x}{\mu\rho}\langle I|[\rho+\frac{x}{(1-x)}](\frac{x}{(1-x)(\rho-x)})+[1+\frac{x}{(1-x)}](\frac{(\rho-\rho x-x)}{(\rho-x)(1-x)})|\rho\rangle$

$=\frac{x}{\mu\rho(1-x)(\rho-x)}\frac{\rho(1-x)(1+x)-x(1-x)}{(1-x)}\langle I|\rho\rangle$

\begin{equation}
=\frac{x}{\mu z\rho(1-x)(\rho-x)}(\rho+\rho x-x).\label{eq:busybusyFinal}\end{equation}

An alternative form:

\begin{equation}
\frac{x}{\mu\rho(1-x)}+\frac{x}{\omega}\label{eq:busybusyAlternative}\end{equation}

can be derived by setting $\hat{g}=\frac{1}{1-xa^{+}}\tilde{g}$ and
then noting that $Q\hat{g}=Q\tilde{g}+x\hat{g}$.

%
{}
\end{document}